\documentclass[10pt]{article}

\usepackage{amsmath}

\usepackage{array}

\usepackage{appendix}

\usepackage{tocloft}                   

\usepackage{graphicx}

\usepackage{amsfonts}

\usepackage{amssymb}

\usepackage{mathrsfs}

\usepackage{yfonts}

\usepackage{euscript}

\usepackage{centernot}                 

\usepackage{ifsym}                     

\usepackage{upgreek}

\usepackage{mathtools}

\usepackage{color}

\usepackage{slantsc}
\usepackage{calligra}

\usepackage{bbold}          

\usepackage[T1]{fontenc}

\usepackage{epsf}

\usepackage{latexsym}

\usepackage{tipa}

\usepackage{makeidx}

\makeindex

\def\be{\begin{equation}}

\def\ee{\end{equation}}

\def\beq{\begin{equation}}

\def\eeq{\end{equation}}

\def\bea{\begin{eqnarray}}

\def\eea{\end{eqnarray}}

\def\ni{\noindent}

\def\!{\hspace{-1.6667em}}

\def\mC{\mbox{C}}                        

\def\mE{\mbox{E}}                        

\def\mI{\mbox{I}}                        

\def\uR{\underline{\mbox{$R$}}}

\def\uq{\underline{\mbox{q}}} 

\def\ur{\underline{\mbox{r}}}

\def\urho{{\underline{\rho}}}

\def\sa{\mbox{\scriptsize a}}

\def\sm{\mbox{\scriptsize m}}

\def\sx{\mbox{\scriptsize x}}

\def\sF{\mbox{\scriptsize F}}

\def\sT{\mbox{\scriptsize T}}

\def\sumi2{\sum\mbox{}_{\mbox{}_{\mbox{\scriptsize $i$=1}}}^2}

\def\sumi3{\sum\mbox{}_{\mbox{}_{\mbox{\scriptsize $i$=1}}}^3}

\def\sumABcycles3{\sum\mbox{}_{\mbox{}_{\mbox{\scriptsize cycles $A,B$=1}}}^{3}}

\def\sumCDcycles3{\sum\mbox{}_{\mbox{}_{\mbox{\scriptsize cycles $C,D$=1}}}^{3}}

\def\sumj3{\sum\mbox{}_{\mbox{}_{\mbox{\scriptsize $j$=1}}}^3}

\def\sumk3{\sum\mbox{}_{\mbox{}_{\mbox{\scriptsize $k$=1}}}^3}

\def\prodiA1{\prod\mbox{}_{\mbox{}_{\mbox{\scriptsize $i$=1}}}^{A - 1}}

\def\d{\textrm{d}}                                                  

\def\FrS{\mbox{\Large $\mathfrak{s}$}}                         

\def\Hilb{\mbox{{\boldmath$\mathfrak{H}$}ilb}}                 
\def\FrQ{\mbox{\Large $\mathfrak{q}$}}                               
\def\Phase{\mbox{{\boldmath$\mathfrak{P}$}hase}}                     

\def\bFrR{\mbox{\boldmath$\mathfrak{R}$}}                            
\def\Rig-Phase{\bFrR\mbox{ig-}\Phase}                                
                                                                													   
\def\FrP{\mbox{\Large $\mathfrak{p}$}}                                 
\def\Positive-Modespace{\mbox{{\boldmath$\mathfrak{M}$}odespace$^+$}}

\def\POSITIVE-MODESPACE{\mbox{{\boldmath$\mathfrak{M}$}ODESPACE$^+$}}

\def\Kin-Hilb{\mbox{{\boldmath$\mathfrak{K}$}in-\Hilb}}                     

\def\Mid-Hilb{\mbox{{\boldmath$\mathfrak{M}$}id-\Hilb}}                     

\def\Dyn-Hilb{\mbox{{\boldmath$\mathfrak{D}$}yn-\Hilb}}                     

\def\5Star{\mbox{\Large$\star$}}              

\def\Frr{\mbox{$\mathfrak{r}$}}

\begin{document}

\begin{titlepage}

\begin{center}

\ni{\bf\Large Maximal Angle Flow on the Shape Sphere of Triangles}

\vspace{.2in}

{\large \bf Edward Anderson$^*$}

\vspace{.1in}

\end{center}

\begin{abstract} 

\ni There has been recent work using Shape Theory to answer the longstanding and conceptually interesting problem of what is the probability that a triangle is obtuse. 
This is resolved by three kissing cap-circles of rightness being realized on the shape sphere; integrating up the interiors of these caps readily yields the answer to be 3/4.  
We now generalize this approach by viewing rightness as a particular value of maximal angle, and then covering the shape sphere with the maximal angle flow. 
Therein, we discover that the kissing cap-circles of rightness constitute a separatrix.
The two qualitatively different regimes of behaviour thus separated both moreover carry distinct analytic pathologies: cusps versus excluded limit points.
The equilateral triangles are centres in this flow, whereas the kissing points themselves -- binary collision shapes -- are more interesting and elaborate critical points. 
The other curves' formulae and associated area integrals are more complicated, and yet remain evaluable.  
As a particular example, we evaluate the probability that a triangle is Fermat-acute, meaning that its Fermat point is nontrivially located; 
the critical maximal angle in this case is 120 degrees.  
 
\end{abstract} 

\vspace{0.1in}

\ni PACS: 04.20.Cv, Physics keywords: configuration spaces, relational spaces, shape spaces, 3-body problem. 

\mbox{ }

\ni Mathematics keywords: Applied Differential Geometry, Applied Topology, Shape Theory, Geometrical Probability, Shape Statistics. 

\vspace{0.1in}
  
\ni $^*$ Dr.E.Anderson.Maths.Physics@protonmail.com

\end{titlepage}

\section{Introduction}

As one of his pillow problems, Charles Lutwidge Dodgson \cite{Pillow} alias Lewis Carroll asked what is the probability that a triangle is obtuse, 
\be
\mbox{ what is } \mbox{ Prob(obtuse) } . 
\ee 
This turned out to be a tough problem -- ambiguous and, moreover, for conceptually interesting reasons -- 
to which many different answers have been given over the years \cite{Guy, Portnoy}.
It has been recently established that \cite{MIT, III, A-Pillow}
\be 
\mbox{Prob(obtuse)} \mbox{ } = \mbox{ } \frac{3}{4} 
\label{3/4}
\ee 
according to Shape Theory in the sense of David Kendall, by use in particular of the shape sphere of triangles 
\cite{Kendall84, Kendall89, Kendall, III} that he discovered.\footnote{See also e.g. 
\cite{Small, GT09, FileR, Bhatta, PE16, KKH16, ABook, I, II, III, IV} for more about Shape Theory.}
%
The underlying interpretation here is that the probability distribution of triangles is the uniform distribution on their configuration space, 
i.e.\ Kendall's shape sphere (as set up in Secs 2 and 3).  
Right triangles moreover trace out three kissing cap-circles thereupon  -- one for each labelling of the hypotenuse -- (Sec 4) 
giving a simple cap area integral from which (\ref{3/4}) follows (Sec 5).  

\mbox{ }

\ni This shape-theoretic answer was first given by Edelman and Strang \cite{MIT}, with a simplified proof in terms of brief and elementary mathematics given in \cite{III, A-Pillow}. 
This is a wide-ranging prototype of mapping flat geometry problems directly realized in a space to shape space, where differential-geometric tools are readily 
available to solve the problem and then finally re-interpret it in the original `shape-in-space' terms.  
Indeed, \cite{A-Pillow} already posed and solved a number of variants of Lewis Carroll's pillow problem. 
Firstly, what is Prob(Isosceles is obtuse).
Secondly, following on from introducing notions of tall and flat triangles so as to investigate the quartet of tall acute, tall obtuse, flat acute and flat obtuse triangles, 
whether isosceles or general.  
The current paper generalizes the scope of pillow problems in a different way, as follows. 

\mbox{ }

\ni We first reconceive of rightness as a critical value of the maximal angle in a triangle, 
\be 
\alpha_{\sm\sa\sx} \mbox{ } = \mbox{ } \frac{\pi}{2} \mbox{ } .
\ee  
This is a commonly encountered critical value due to qualitative differences between obtuse triangles 
\be
\alpha_{\sm\sa\sx} \mbox{ } >  \mbox{ } \frac{\pi}{2} \mbox{ } , 
\ee
and acute triangles
\be
\alpha_{\sm\sa\sx} \mbox{ } <  \mbox{ } \frac{\pi}{2} \mbox{ } , 
\ee
many simplifications occur for right angles as well.\footnote{For instance, the cosine rule collapses to Pythagoras' Theoerem, 
and the chord subtending the angle on a circle is now the diameter.}
%
It it not however the only critical value of $\alpha_{\sm\sa\sx}$ as regards qualitative change in geometrical properties.

\mbox{ }

\ni Let us first recollect that the allowed range of $\alpha_{\sm\sa\sx}$ is  
\be
\frac{\pi}{3} \mbox{ } \leq \mbox{ } \alpha_{\sm\sa\sx} \mbox{ } \leq \mbox{ } \pi \mbox{ } . 
\ee 
This follows from the angle sum of a triangle being $\pi$. 
The strict case of the first inequality is of course the equilateral triangle, 
whereas the second returns a degenerate case (excluded in some applications): the collinear triangles without coincident vertices . 

\mbox{ }

\ni So for instance, 
\be 
\alpha_{\sm\sa\sx} \mbox{ } = \mbox{ } \frac{2 \, \pi}{3} 
\ee 
is a meaningful critical value as regards Fermat's problem (Fig \ref{Fermat-Points}.a)\footnote{Whereas Fermat posed this problem, it was Torricelli who solved it.  
So this point is sometimes known as the Torricelli or Fermat--Torricelli point.  
See \cite{17th} for these 17th century origins of this problem.} 
of which point $F$ minimizes the total distance from itself to a triangle's three vertices (Fig \ref{Fermat-Points}).  
This is a critical value for Fermat's problem because for 
\be
\alpha_{\sm\sa\sx} \mbox{ } \geq  \mbox{ } \frac{2 \, \pi}{3} \mbox{ } , 
\ee 
the Fermat point $F$ is just at the obtuse angle's vertex (Fig \ref{Fermat-Points}.b), whereas for 
\be
\alpha_{\sm\sa\sx} \mbox{ } < \mbox{ } \frac{2 \, \pi}{3} \mbox{ } , 
\ee 
the Fermat point $F$ is a nontrivial point in the interior of the triangle (Fig \ref{Fermat-Points}.c).  
These are clearly qualitatively different regimes.
On these grounds, we make the following definitions.
%
{            \begin{figure}[!ht]
\centering
\includegraphics[width=0.7\textwidth]{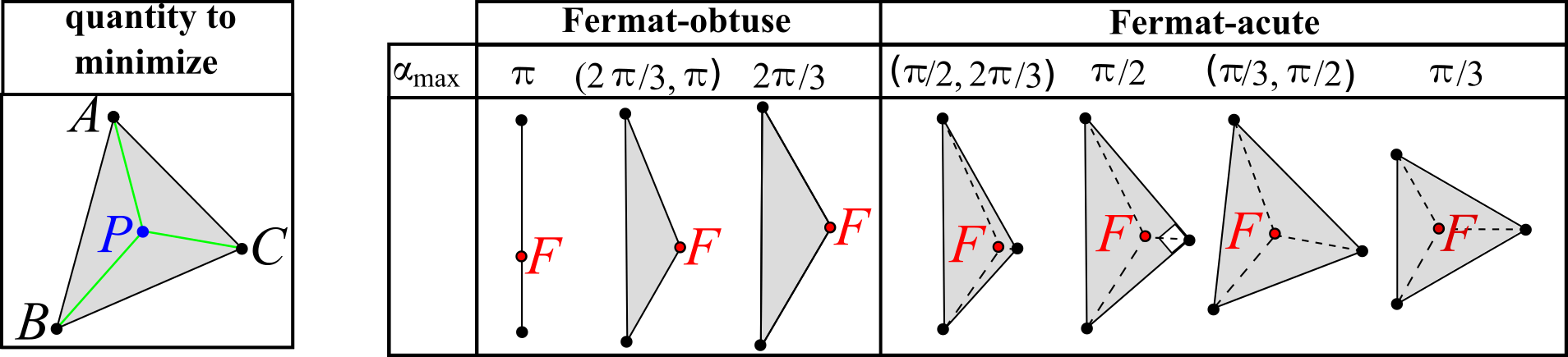}
\caption[Text der im Bilderverzeichnis auftaucht]{        \footnotesize{a) Fermat's problem is, given a triangle ABC, to find the position of the point P (in blue) 
such that the sum of the lengths (in green) from it to the three vertices is minimized. 
b) This is solved by the Fermat point $F$ (in red), which is moreover just at the obtuse vertex for triangles with $\alpha_{\sm\sa\sx} \geq 2 \, \pi/3$, which we term Fermat-obtuse. 
On the other hand, the Fermat point $F$ has a nontrivial position for triangles with $\alpha_{\sm\sa\sx} <  2 \, \pi/3$.  
In this case, $F$ is such that the vertices are spread out at $2 \, \pi/3$ from each other relative to $F$, as indicated by the dashed lines. 
Moreover, for $\alpha_{\sm\sa\sx} =  \pi/3$, the position of $F$ reduces to just that of the centre of symmetry of the equilateral triangle.} }
\label{Fermat-Points}\end{figure}            }

\mbox{ }  

\ni{\bf Definition 1} A triangle is 
 \be 
\mbox{{\it Fermat-acute} if } \mbox{ }   \mbox{ } \alpha_{\sm\sa\sx}   \mbox{ }  <  \mbox{ } \frac{2 \, \pi}{3} \mbox{ } ,
\ee
\be 
\mbox{{\it Fermat-critical} if } \mbox{ } \mbox{ } \alpha_{\sm\sa\sx} = \frac{2 \, \pi}{3} \mbox{ } ,  \mbox{ } \mbox{ and } 
\ee
\be 
\mbox{{\it Fermat-obtuse} if }  \mbox{ }  \mbox{ } \alpha_{\sm\sa\sx}   \mbox{ } >  \mbox{ } \frac{2 \, \pi}{3} \mbox{ } . 
\ee
\ni{\bf Definition 2} A triangle is 
 \be 
\mbox{{\it $\alpha$-acute} if }    \mbox{ }  \mbox{ }  \alpha_{\sm\sa\sx} \mbox{ } < \alpha \mbox{ } , 
\ee
\be 
\mbox{{\it $\alpha$-critical} if } \mbox{ }  \mbox{ }  \alpha_{\sm\sa\sx} = \alpha \mbox{ } , \mbox{ } \mbox{ and } 
\ee
\be 
\mbox{{\it $\alpha$-obtuse} if }   \mbox{ }  \mbox{ }  \alpha_{\sm\sa\sx}  \mbox{ } > \mbox{ } \alpha \mbox{ } . 
\ee
This case covers both Fermat-criticality and right-criticality, 
while being an open invitation for others to point to values of $\alpha_{\sm\sa\sx}$ that are critical in further geometrical contexts, 
by which to extend the current paper's range of examples.    

\mbox{ }

\ni In this manner, the following further pillow problems can be set. 
Firstly, 
\be 
\mbox{what is } \mbox{ Prob(Fermat-obtuse)}?  
\ee 
is the most direct counterpart of the original pillow problem, though, 
since it is its complement being the geometrically nontrivial case, one might primarily pose this instead: 
\be 
\mbox{what is } \mbox{ Prob(Fermat-acute)}?  
\ee
Secondly, 
\be 
\mbox{what is } \mbox{ Prob($\alpha$-obtuse)}?  
\label{A-obtuse}
\ee
\ni In the current paper, we provide a solution for this problem by finding the maximal angle flow on the shape sphere (Sec 7).
This extends the kissing caps of right-maximality to a congruence of curves covering all other maximal angles. 
We show that for rightness the calculation exceptionally collapses, 
yet the other curves remain calculable, as are the corresponding enclosed areas evaluating (\ref{A-obtuse}), given in Sec 8 for the particular Fermat example.  

\mbox{ } 

\ni We show moreover in Sec 7 that the kissing cap-circles of rightness plays the furtherly critical role of a separatrix in this maximal angle flow. 
The two qualitatively different regimes of behaviour this separates both moreover carry distinct analytic pathologies: cusps versus excluded limit points.
The flow also contains three copies of an unusual type of critical point at the kissing points themselves. 
Thereby, this flow is likely to be of interest in Differential Geometry and Dynamical Systems, as well as through its meaningfulness for Flat Geometry.  
We finally comment on further Flat Geometry problems opened up by Kendall's Shape Theory in the Conclusion (Sec 9).

\vspace{10in}

\section{A progression of coordinatizations for the general triangle}
%
{            \begin{figure}[!ht]
\centering
\includegraphics[width=1.0\textwidth]{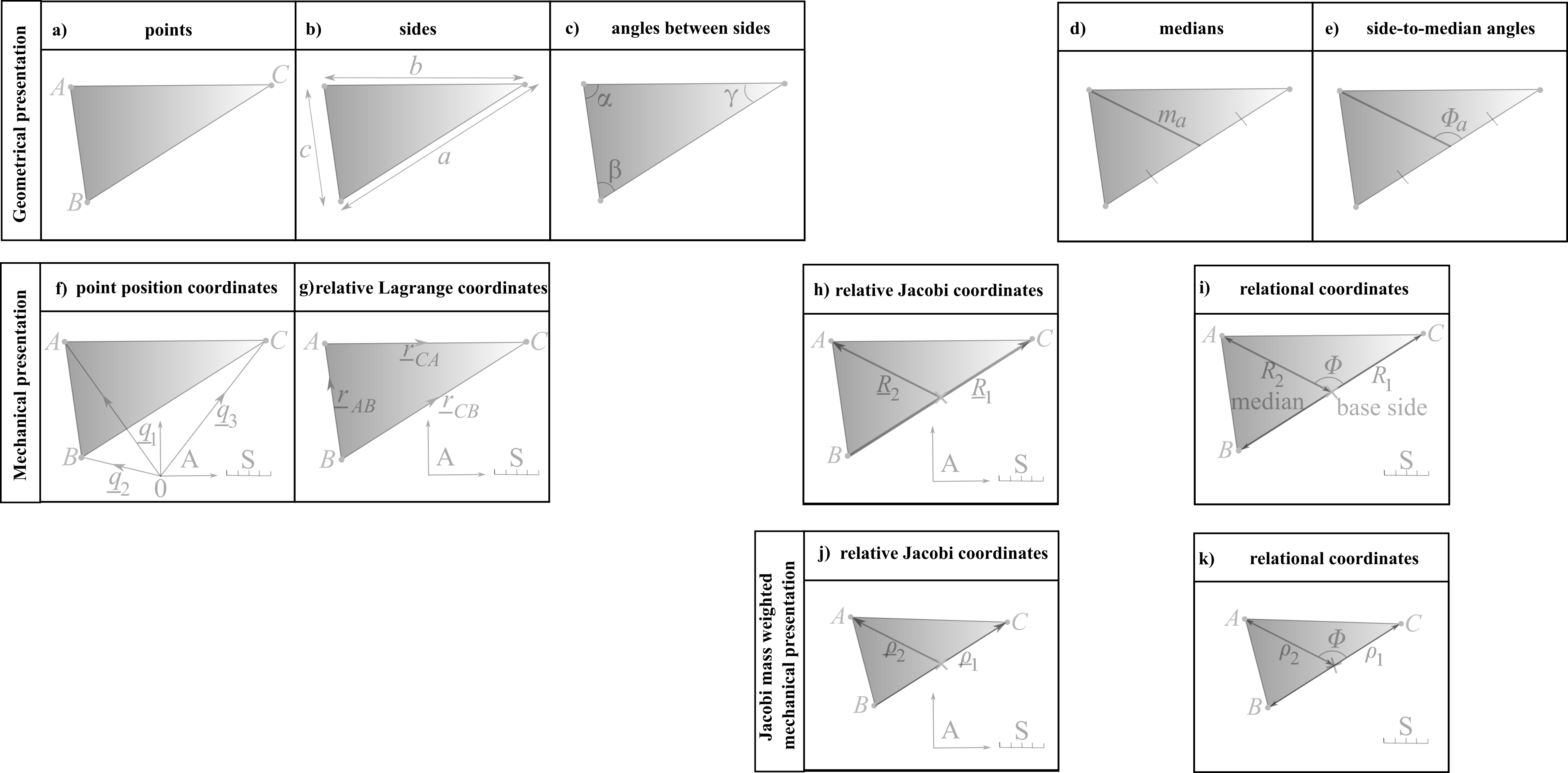}
\caption[Text der im Bilderverzeichnis auftaucht]{        \footnotesize{Geometrical and mechanical notation for the general triangle.} }
\label{Tri-q-r-rho-relatio}\end{figure}            }
%
\ni Let us first rephrase Fig \ref{Tri-q-r-rho-relatio}.a)-c)'s standard geometrical description of the general triangle 
in terms of Fig \ref{Tri-q-r-rho-relatio}.f)'s {\it position coordinates} $\uq_I$ for the triangle's vertices -- now viewed as particles -- 
relative to an absolute origin 0, axes A and scale S.

\mbox{ }

\ni Let us next pass to {\it relative Lagrangian coordinates} (Fig \ref{Tri-q-r-rho-relatio}.e), consisting of particle separation vectors 
\be
\ur_{IJ} = \uq_J - \uq_I \mbox{ } .  
\ee 
Using these amounts to discarding the absolute origin 0.  
Their magnitudes return the side-lengths and the angles between these vectots return Fig \ref{Tri-q-r-rho-relatio}.c)'s angles.  

\mbox{ }

\ni Moreover, not all of the relative Lagrange coordinates are independent. 
This can be circumvented by picking a basis of two of them.  
This however leaves the inertia quadric in non-diagonal form.  

\mbox{ }

\ni Diagonalizing the inertia quadric leads us on to {\it relative Jacobi coordinates} \cite{Marchal} (Fig \ref{Tri-q-r-rho-relatio}.h).
These are a more general notion of particle {\sl cluster} separation vectors, i.e.\  particle subsystem centre of mass separations.   
They include the relative Lagrange coordinates as a special subcase: the one for which the particle subsystems at both ends consist of one particle each, 
by viewing this particle as the location of its own centre of mass.
In particular, in the case of a triangular configuration, 
\be
\uR_1          =          \uq_C - \uq_B                   \mbox{ } , \mbox{ } \mbox{ } 
\uR_2 \mbox{ } = \mbox{ } \uq_A - \frac{\uq_B + \uq_C}{2} \mbox{ }
\label{R-def}
\ee 
or cycles thereof for the other two choices of 2-particle cluster.  
These have associated cluster masses (conceptually reduced masses) 
\be 
\mu_1 \mbox{ } = \mbox{ } \frac{1}{2} \mbox{ } , \mbox{ } \mbox{ } \mu_2 \mbox{ } = \mbox{ } \frac{2}{3}
\ee 
respectively. 
For later use, note furthermore that 
\be
\sqrt{\frac{\mu_2}{\mu_1}} = \sqrt{\frac{\frac{2}{3}}{\frac{1}{2}}} = \frac{2}{\sqrt{3}} \mbox{ } . 
\label{mu-ratio}
\ee
\ni One can additionally pass to mass-weighted relative Jacobi coordinates (Fig \ref{Tri-q-r-rho-relatio}.j)
\be 
\urho_i := \sqrt{\mu_i}\uR_i \mbox{ } \mbox{ } (i = 1 \mbox{ } , \mbox{ } \mbox{ } 2) \mbox{ } ; 
\label{rho-R}
\ee 
These leave the inertia quadric's matrix as the identity, and furthemore turn out to be shape-theoretically convenient to work with.

\mbox{ }

\ni Passing to mass-weighted relative Jacobi coordinates $\urho_a$ does not moreover do anything as regards removing the absolute axes A or scale S, 
since these remain vectorial, and thus made reference to absolute directions and absolute magnitudes. 
Reference to absolute axes A is moreover removed by restricting attention [Fig \ref{Tri-q-r-rho-relatio}.i) and k)] to the following triple.

\mbox{ }

\ni i) The {\it relative Jacobi magnitudes} $\rho_a$ (a = 1, 2); comparing Fig \ref{Tri-q-r-rho-relatio}.a) and d) with the mass-unweighted version of Fig \ref{Tri-q-r-rho-relatio}.h), 
these are a base side and the corresponding median respectively.
 
\mbox{ } 
 
\ni ii) The {\it relative Jacobi angle} between the $\urho_a$, 
\be 
\Phi      \mbox{ } = \mbox{ } \mbox{arccos}    \left(    \frac{\urho_1\cdot\urho_2}{\rho_1\rho_2}  \right)
          \mbox{ } = \mbox{ } \mbox{arccos}    \left(    \frac{\sqrt{\mu_1} \uR_1 \cdot \sqrt{\mu_2} \uR_2}{ \sqrt{\mu_1} R_1 \sqrt{\mu_2} R_2} \right) 
          \mbox{ } = \mbox{ } \mbox{arccos}    \left(    \frac{\uR_1 \cdot \uR_2}{ R_1 R_2}     \right) \mbox{ } . 
\ee  
This angle can be thought of more vividly as a `{\it Swiss-army-knife' relative angle} \cite{FileR}.    
Note that the above demonstrates this angle to be unaltered by mass-weighting, 
by which it is identical to the geometrical side-to-median angle $\Phi_a$ of Fig \ref{Tri-q-r-rho-relatio}.e)

\mbox{ }

\ni Together, $\rho_1, \rho_2, \Phi$ constitute {\it relational scale-and-shape data} for the triangle \cite{III}.

\mbox{ }

\ni Finally, reference to absolute scale S is removed by \cite{FileR} continuing to consider the relative Jacobi angle, but now just alongside 
the {\it ratio of Jacobi magnitudes},  
\be 
{\cal R} \mbox{ } := \mbox{ } \frac{\rho_2}{\rho_1} \mbox{ } .
\label{cal-R-def}
\ee 
In \cite{III} we also ascertained that this provides a further similarity condition for triangles and, accordingly, serves as data from which 
the `primary' triangle's features of Fig 1.a) modulo scale -- side-to-side relative angles and side-length ratios -- can be abstracted.

\section{Kendall's shape sphere of triangles}
%
{            \begin{figure}[!ht]
\centering
\includegraphics[width=1.0\textwidth]{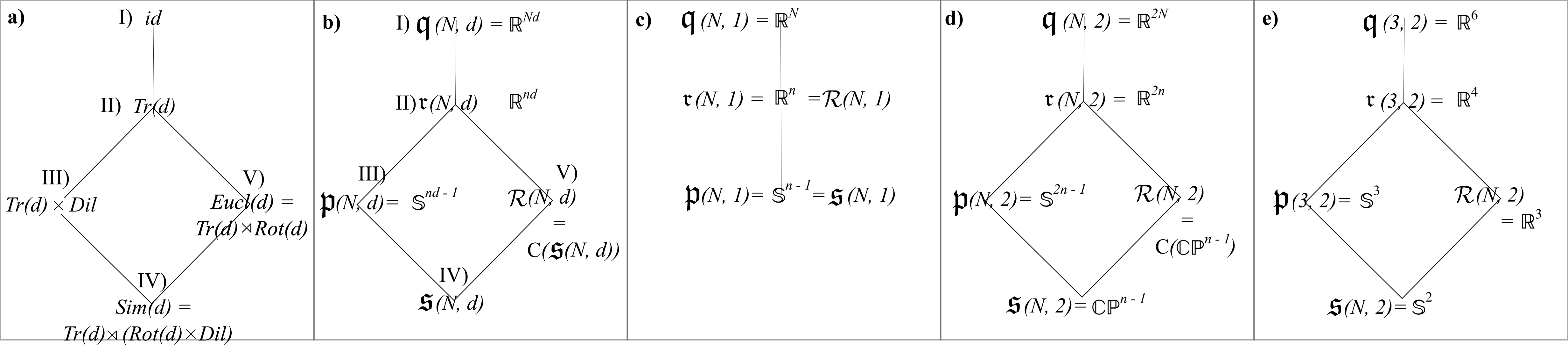}
\caption[Text der im Bilderverzeichnis auftaucht]{        \footnotesize{a) Lattice of 5 subgroups of the similarity group $Sim(d)$, 
with b) the corresponding lattice of configuration space quotients, specializing to dimension 1 in c), 2 in d) and furthermore to the 3 particles in 2-$d$ of triangleland in e).} }
\label{Lattice-of-5}\end{figure}            }

\mbox{ }

\ni I) At the level of configuration spaces, the space of the $\uq_I$ for arbitrary particle number $N$ and dimension $d$ is straightforwardly 
\be 
\FrQ(N, d) = \{\mathbb{R}^d\}^N = \mathbb{R}^{N \, d}                      \mbox{ } . 
\ee 
Let us next follow through what happens to this as one removes the absolute origin 0, axes A and scale S \cite{Kendall84, Kendall, FileR}. 

\mbox{ }

\ni II) Quotienting out by $Tr(d)$ -- the $d$-dimensional translations --  we arrive at the configuration space of independent $\ur_{IJ}$, $\uR_i$ or $\urho_i$: {\it relative space} 
\be 
\Frr(N, d) \mbox{ } := \mbox{ } \frac{\FrQ(N, d)}{Tr(d)} 
           \mbox{ }  = \mbox{ } \frac{\mathbb{R}^{N \, d}}{\mathbb{R}^d} 
		   \mbox{ }  = \mbox{ } \mathbb{R}^{n \, d}                       \mbox{ } , 
\ee
where we have also made use of $n := N - 1$.  

\mbox{ }

\ni III) It is next useful to consider quotienting out the dilation group $Dil$ to render absolute scale S irrelevant. 
This yields Kendall's {\it preshape space} \cite{Kendall84, Kendall} 
\be 
\FrP(N, d) \mbox{ } := \mbox{ } \frac{\Frr(N, d)}{Dil} 
           \mbox{ }  = \mbox{ } \frac{\mathbb{R}^{n \, d}}{\mathbb{R}_+} 
           \mbox{ }  = \mbox{ } \mathbb{S}^{n \, d - 1}                   \mbox{ } :   		   
\ee
the ($n \, d - 1$)-sphere, which is the configuration space of ratios of independent components of $\ur_{IJ}$, $\uR_i$ or $\urho_i$.
This moreover carries the standard hyperspherical metric.

\mbox{ }

\ni IV) Quotienting out as well by $Rot(d)$ -- the $d$-dimensional rotations -- 
and so overall by $Sim(d)$: the $d$-dimensional similarity group of translations, rotations and dilations. 
we arrive at Kendall's {\it shape space} \cite{Kendall84, Kendall89, Kendall} is 
\be 
\FrS(N, d) \mbox{ } := \mbox{ } \frac{\FrQ(N, d)}{Sim(d)} 
           \mbox{ }  = \mbox{ } \frac{\Frr(N, d)}{Rot(d) \times Dil} 
		   \mbox{ }  = \mbox{ } \frac{\FrP(N, d)}{Rot(d)}                 \mbox{ } . 
\ee 
\ni Subsequent analysis picks up $d$ and $N$ dependence.
For $d = 1$, there are not continuous rotations to remove, so 
\be 
\FrS(N, 1) = \FrP(N, 1) = \mathbb{S}^{n - 1}                             \mbox{ } . 
\ee 
For $d = 2$, 
\be
Rot(2) = SO(2) = U(1) = \mathbb{S}^1 \mbox{ } , 
\ee 
and \cite{Kendall84} 
\be 
\FrS(N, 2) \mbox{ } = \frac{ \mathbb{S}^{2 \, n - 1}}{\mathbb{S}^1} 
           \mbox{ } = \mbox{ } \mathbb{CP}^{n - 1}                       \mbox{ } , 
\ee 
which is most readily arrived at by the complex Hopf map generalization of the Hopf map \cite{Frankel}.  
Moreover, shape space carries the natural Fubini--Study metric \cite{Kendall84, Kendall}. 
This result more than suffices for the current paper; indeed, for $N = 3$ triangleland, the Hopf map 
\be
\mathbb{S}^3 \longrightarrow \mathbb{S}^2 \mbox{ } ( \mbox{ } = \mbox{ } \mathbb{CP}^1 \mbox{ } ) \mbox{ } 
\ee 
suffices, and the Fubini--Study metric simplifies in plane-polar coordinates to \cite{Kendall, FileR}
\be 
\d s^2 = 4 \frac{\d {\cal R}^2 + {\cal R}^2\d\Phi^2}{  (  1 + {\cal R}^2  )^2  } \mbox{ } .   
\ee 
This can in turn be recognized as the standard spherical metric in stereographic coordinates. 
Thus, in the Shape Theory of triangles, the ratio of Jacobi magnitudes ${\cal R}$ of eq. (\ref{cal-R-def}) plays the geometric role of stereographic radius of the shape sphere. 

\mbox{ }
 
\ni The Swiss-army-knife relative angle $\Phi$ of the shape-in-space thus plays the role of polar angle on the shape sphere. 

\mbox{ }

\ni Finally, the venerable substitution 
\be 
{\cal R} = \mbox{tan} \, \frac{\Theta}{2} \mbox{ }  
\label{RTheta}
\ee 
serves to convert the stereographic radius to the standard azimuthal coordinate $\Theta$ (about the U-axis), casting the shape sphere metric into the standard spherical metric form 
\be 
\d s^2 = \d {\Theta}^2 + \mbox{sin}^2 \Theta \, \d\Phi^2 \mbox{ } . 
\ee
\ni See Fig \ref{S(3, 2)-Intro} for where some of the most qualitatively distinctive triangle shapes-in-space are realized in the triangleland shape space.
%
{            \begin{figure}[!ht]
\centering
\includegraphics[width=0.45\textwidth]{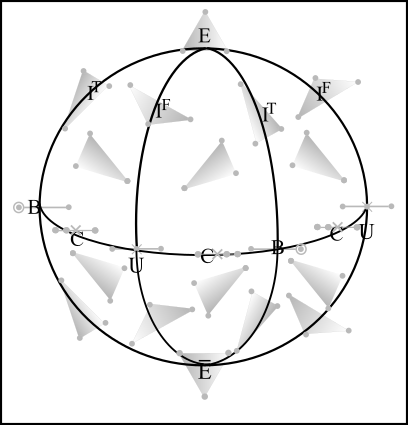}
\caption[Text der im Bilderverzeichnis auftaucht]{        \footnotesize{a) The triangleland shape sphere \cite{Kendall89, +Tri, FileR, ABook, III}.  
Equilateral triangles E are at its poles, whereas collinear configurations C form its equator.
There are 3 bimeridians of isoscelesness I corresponding to 3 labelling choices for the vertices.
For later use in this paper, I make further distinction between the tall isosceles triangles $\mI^{\sT}$ and the flat isosceles triangles $\mI^{\sF}$.  
$\mC \, \cap \mI$ gives 3 binary collisions B \cite{Montgomery}, which are moreover topologically significant, and 3 uniform collinear shapes U.
Each of these triples lies on the equator, with each member at $2 \, \pi/3$, and the two triples at $\pi/3$ to each other.} }
\label{S(3, 2)-Intro} \end{figure}          }

\mbox{ }

\ni Note that the spherical coordinates we have found have U and B as their North and South poles.
There are moreover 3 clustering choices of such axes, realized at  $2 \, \pi /3$ angles to each other in the plane of collinearity.  
For some purposes, one would prefer to use the more distinguished and cluster-independent $\mE\overline{\mE}$ axis to define new spherical polar coordinates 
$\widetilde{\Theta}$, $\widetilde{\Phi}$ \cite{FileR, III}.  
Working in terms of these is however more involved in terms of moving back and forth between shapes in space and points, curves and regions in the shape sphere.  
Finally note that we only need one set of spherical polars to establish the shapes-in-space to points in shape space correspondence.  

\vspace{10in}

\ni V) Finally, quotienting out instead by $Eucl(d)$ -- the Euclidean group of translations and rotations -- 
we arrive at the topological and geometrical form of the {\it relational space} \cite{FileR, I}
\be 
{\cal R}(N, d) \mbox{ } := \mbox{ } \frac{\FrQ(N, d)}{Eucl(d)} 
               \mbox{ }  = \mbox{ } \frac{\Frr(N, d)}{Rot(d)}     \mbox{ } . 
\ee 
\ni A general result for this is that 
\be 
{\cal R}(N, d) = \mC(\FrS(N, d)) \mbox{ } , 
\ee 
for C the topological and metric coning construct.
In particular, for the 2-$d$, $N = 3$ triangleland, 
\be 
{\cal R}(3, 2) = \mC(\mathbb{S}^2) = \mathbb{R}^3
\ee 
topologically (albeit it is not metrically flat, see e.g. \cite{FileR}).

\vspace{10in}

\vspace{10in}

\section{The kissing cap-circles of rightness}

\ni We now address the question of where in the shape sphere the right triangles are located.

\mbox{ }

\ni{\bf Theorem 1} \cite{MIT, III, A-Pillow} The right triangles form three equal cap-spheres, each with azimuthal coordinate $\frac{\pi}{3}$ about the corresponding U point. 
These kiss at each B point.  

\mbox{ }

\ni{\underline{Proof}} 
Begin with the very well-known Theorem of Euclidean Geometry that the angle subtended by a diameter of a circle is right.
In the case of a right triangle's circumcircle (Fig \ref{2-to-1}.a), the hypotenuse is a diameter and its corresponding median is a radius. 
Thus, hypotenuse : median = diameter : radius = 2 : 1.
So, in relative Jacobi magnitudes as per the identification under (\ref{R-def}), 
\be 
\frac{R_2}{R_1} \mbox{ } = \mbox{ } \frac{1}{2} \mbox{ } .
\label{half}
\ee 
%
{            \begin{figure}[!ht]
\centering
\includegraphics[width=0.45\textwidth]{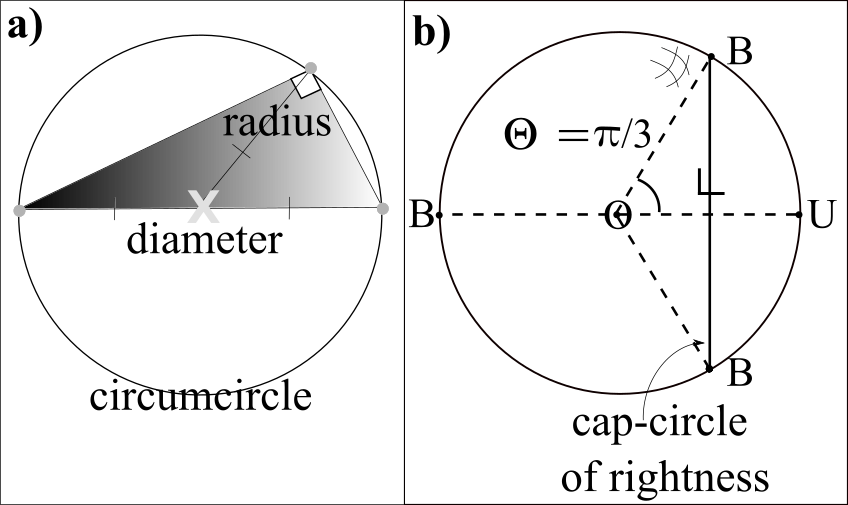}
\caption[Text der im Bilderverzeichnis auftaucht]{        \footnotesize{a) Right-angled triangles' hypotenuses and their corresponding medians are in 2 : 1 ratio. 

\mbox{ }

\ni b) Each labelling of right triangles correspond to a cap-circle of constant azimuthal value $\frac{\pi}{3}$.} }
\label{2-to-1} \end{figure}          }

\mbox{ }

\ni Next insert mass-weighting to obtain 
\be
{\cal R} \mbox{ } = \mbox{ } \frac{\rho_2}{\rho_1} 
         \mbox{ } = \mbox{ } \sqrt{\frac{\mu_2}{\mu_1}} \, \frac{R_2}{R_1}
		 \mbox{ } = \mbox{ } \frac{2}{\sqrt{3}} \, \frac{1}{2}
		 \mbox{ } = \mbox{ } \frac{1}{\sqrt{3}} \mbox{ } .  
		 \label{R=1/sqrt3}
\ee
Here we have used (\ref{cal-R-def})                  in the first equality, 
                   the magnitudes of (\ref{rho-R})   in the second, 
				   (\ref{mu-ratio}) and (\ref{half}) in the third, 
				   and an elementary cancellation    in the fourth. 

\mbox{ }

\ni Furthermore, in spherical polar coordinates, 
\be 
\mbox{tan}\frac{\Theta}{2}  \mbox{ } = \mbox{ } {\cal R} 
							\mbox{ } = \mbox{ } \frac{1}{\sqrt{3}}      \mbox{ } .
\ee
Thus, inverting,  
\be
\Theta \mbox{ } = \mbox{ } 2 \times \frac{\pi}{6} 
       \mbox{ } = \mbox{ } \frac{\pi}{3}                                \mbox{ } . 
\label{Theta=pi/3}
\ee  
This is sketched in Fig \ref{2-to-1}.b) for one cluster choice's cap-circle of rightness. 
Repeat for each of the three clusters to obtain three cap-circles. 
Each is centred about the corresponding U-point pole, and kisses the other two at the two B points which are not antipodal to that U-point: 
{\it kissing cap-circles of rightness} as viewed in Fig \ref{Kissing} from a variety of directions.                                            $\Box$

\mbox{ }

\ni The regions outside of these caps are acute, seen e.g.\ from the observation that the equilateral triangles E residing at the poles are outside of these caps and are acute, 
and then applying a continuity argument.  
Conversely, the interiors of the caps are the obtuse triangles, so we term the caps themselves the {\it kissing caps of obtuseness}. 
There being three caps is commeasurate with there being three ways a labelled triangle can be obtuse.
%
{            \begin{figure}[!ht]
\centering
\includegraphics[width=0.7\textwidth]{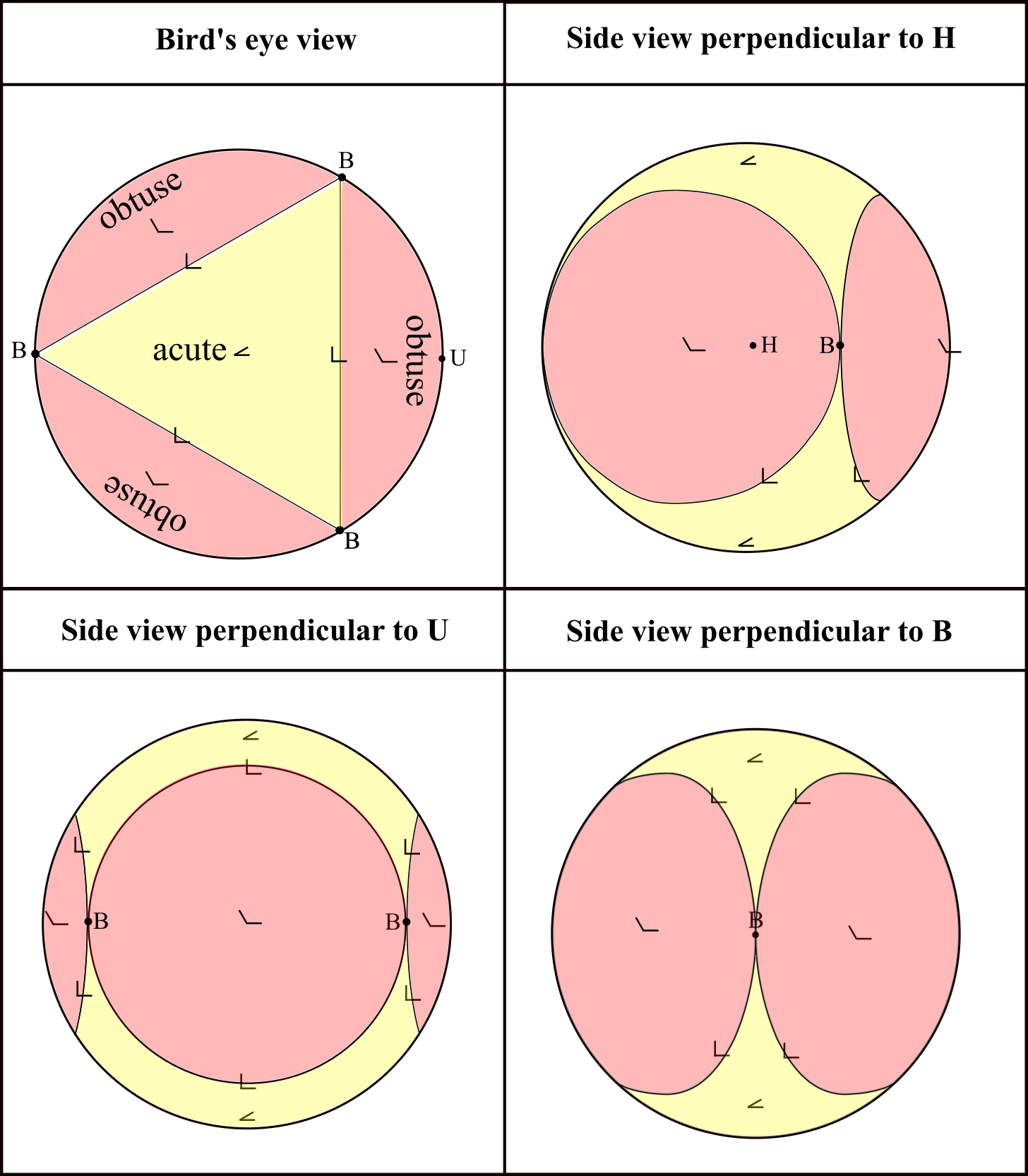}
\caption[Text der im Bilderverzeichnis auftaucht]{        \footnotesize{Kissing cap-circles of rightness, bounding kissing caps of obtuseness, as viwed from various directions.
H is the collinear shape `halfway' between B and U \cite{III}; it lies in the direction perpendicular to the UB axis. } }
\label{Kissing} \end{figure}          }

\mbox{ }

\mbox{ }

\mbox{ }

\mbox{ }

\mbox{ }

\mbox{ }

\mbox{ }

\mbox{ }

\mbox{ }

\mbox{ }

\mbox{ }

\mbox{ }

\mbox{ }

\mbox{ }

\mbox{ }

\mbox{ }

\vspace{10in}

\section{Shape-theoretic solution of Lewis Carroll's pillow problem}

{\bf Corollary 1} i) 
\be
\mbox{Prob}(\mbox{obtuse}) = \mbox{ } \frac{3}{4}   \mbox{ } . 
\ee
ii)
\be
\mbox{Prob}(\mbox{acute}) = \mbox{ } \frac{1}{4}   \mbox{ } . 
\ee
\ni{\underline{Proof}} 
\ni Each cap of obtuseness contributes an area (or more basically, surface of revolution) integral\footnote{Archimedes 
gave an even more basic derivation for the area of a spherical cap, by which Calculus could be avoided altogether in this derivation.  } 
\be
\int_{\Phi = 0}^{2\,\pi} \int_{\Theta = 0}^{\frac{\pi}{3}} \mbox{sin} \, \Theta \, \d \Theta \, \d \Phi  \mbox{ } = \mbox{ } 
2 \, \pi \int_{\Theta = 0}^{\frac{\pi}{3}} \mbox{sin} \, \Theta \, \d \Theta 
                                                                                                          = 2 \, \pi \left[ - \mbox{cos} \, \Theta \right]_{\Theta = 0}^{\frac{\pi}{3}}  
																				                 \mbox{ } = \mbox{ } 2 \, \pi \left( - \frac{1}{2} + 1 \right)                                                           
																								 \mbox{ } = \mbox{ }      \pi                                             \mbox{ } . 
\ee
Thus, between them, the three caps of obtuseness contribute an area of $3 \, \pi$.
On the other hand, the area of the whole 2-sphere is $4 \, \pi$, so 
\be
\mbox{Prob}(\mbox{obtuse}) = \frac{3 \, \pi}{4 \, \pi} \mbox{ } = \mbox{ } \frac{3}{4}   \mbox{ } . 
\ee
\ni Finally, complementarily 
\be
\mbox{Prob}(\mbox{acute})          =           1 - \mbox{Prob}(\mbox{obtuse})         
                          \mbox{ } = \mbox{ }  1 - \frac{3}{4}               
						  \mbox{ } = \mbox{ }      \frac{1}{4}                           \mbox{ } . \mbox{ } \Box 
\ee
\ni Note that this answer is in the context of interpreting the question's probability distribution to be the {\sl uniform} one {\sl on the corresponding shape space}, 
as equipped with the standard spherical metric that Kendall \cite{Kendall84, Kendall89, Kendall} showed to be natural thereupon. 
The uniqueness of this metric is furtherly intuitively obvious from its being induced as a quotient of a simpler structure (position space or relative space: flat spaces), 
or from its arising correspondingly by reduction of the corresponding mechanical actions. 

\mbox{ }
 
\ni Also note that 3/4 is in fact one \cite{Portnoy} of the more common answers to Lewis Carroll's pillow problem, now obtained moreover on shape-theoretic premises. 
Edelman and Strang \cite{MIT} recently obtained this result using Kendall's Shape Theory, 
while \cite{A-Pillow} first gave a simplified proof along the lines provided above. 
In conjunction with various elementary derivations of Kendall's shape sphere \cite{A-Pillow, Forthcoming}, this reduces the shape-theoretic answer to Lewis Carroll's pillow 
problem to just elementary algebra and (pre)calculus.  
This renders it open,     both to anybody studying STEM subjects at University (or even High School in many countries) 
                      and also to considerable generalization to a whole class of geometrical problems, 
					  of which \cite{A-Pillow}, the current paper and \cite{IV} provide the first few. 

\mbox{ }

\ni On the one hand, see also \cite{Guy} for various other answers' values for $\mbox{Prob(obtuse)}$. 

\mbox{ }

\ni On the other hand, \cite{Portnoy} moreover asked for a general principle behind the methods which share the answer 3/4. 
Shape Theory as per above surely provides a principle,  
though whether {\sl all} known methods which yield 3/4 can be shown to follow from this shape-theoretic principle is left as a good topic for a further paper.  

\vspace{10in}

\section{The maximum angle flow}

The current paper's first new result is that the maximal angle curves on the shape sphere are underlied by the following constant-angle curves.  

\mbox{ }

\ni{\bf Theorem 2} The curves of constant angle 
\be
\alpha \neq \frac{\pi}{2}
\ee 
-- the angle defined in Fig \ref{Tri-q-r-rho-relatio}.a) -- are given by  
\be
\Phi   \mbox{ } = \mbox{ }   \mbox{arcsin} \left(  \frac{    1 - 2 \, \mbox{cos} \, \Theta    }{    \sqrt{3} \, k \, \mbox{sin} \, \Theta }  \right) \mbox{ } ,  
\label{Phi-eq-2}
\ee 
for 
\be 
k := \mbox{cot} \, \alpha \mbox{ } . 
\ee 
\ni{\underline{Proof}}
The rescaling 
\be 
{\cal V} = \sqrt{3} \, {\cal R}
\label{VR}
\ee  
is useful \cite{III} in solving for triangles in terms of shape space data.
In terms of this, \cite{III} establishes with a sequance of basic trigonometric moves that 
\be 
\mbox{cot} \, \alpha \mbox{ } = \mbox{ } \frac{{\cal V}^2 - 1}{2 \, {\cal V} \, \mbox{sin} \, \Phi} \mbox{ } .  
\ee 
We now make $\Phi$ the subject of this as follows. 
Firstly, 
\be
2 \, k \, \mbox{sin} \, \Phi = {\cal V} - {\cal V}^{-1} \mbox{ } . 
\ee
There are then two mutually-exhaustive cases to consider. 

\mbox{ }

\ni Case 1) For $k = 0$, i.e.\ $\alpha = \pi/2$: right-angled triangles, we have 
\be
{\cal V}^2 = 1 \mbox{ } . 
\label{k=0}
\ee
Since ${\cal V}$'s definition implies non-negativity, this means that 
\be 
{\cal V} =  1 \mbox{ } \mbox{ } \Rightarrow \mbox{ } \mbox{ } {\cal R} = \frac{1}{\sqrt{3}}  \mbox{ } , 
\ee 
amounting to a recovery of Sec 4's result (\ref{R=1/sqrt3}). 

\mbox{ }

\ni Case 2) For $k \neq 0$, i.e.\ $\alpha = \pi/2$, 
\be
\Phi = \mbox{arcsin}\left(\frac{1}{2 \, k}\left({\cal V} - {\cal V}^{-1}\right)\right)            \mbox{ } .  
\label{Phi-eq}
\ee 
\ni Then 
\be 
{\cal V} - {\cal V}^{-1} \mbox{ } = \mbox{ } \sqrt{3} \, \mbox{tan} \, \frac{\Theta}{2} - \frac{1}{\sqrt{3}}  \mbox{cot} \, \frac{\Theta}{2} 
                         \mbox{ } = \mbox{ } \frac{          3  \, \mbox{sin}^2  \frac{\Theta}{2} -  \mbox{cos}^2 \frac{\Theta}{2}    }
						              {    \sqrt{3} \, \mbox{sin} \, \frac{\Theta}{2} \, \mbox{sin} \, \frac{\Theta}{2}     } 
						 \mbox{ } = \mbox{ }       \frac{    1 - 2 \, \mbox{cos} \, \Theta    }{    \frac{\sqrt{3}}{2} \, \mbox{sin} \, \Theta }			 \mbox{ } .   
\ee 
Finally substitute this in (\ref{Phi-eq}) to obtain the Theorem. $\Box$

\mbox{ }

\ni{\bf Corollary 2} Introducing the shape sphere's counterpart of the {\it Legendre variable}, 
\be 
{\cal X} := \mbox{cos} \, \Theta \mbox{ } ,  
\ee 
the further simplification   
\ni 
\be
\Phi      = \mbox{arcsin}  \left(  \frac{    1 - 2 \, {\cal X}    }{    \sqrt{3} \, k \, \sqrt{   1 - {\cal X}^2  }   }  \right) 
\label{Phi-eq-3}
\ee 
is afforded.
%
{            \begin{figure}[!ht]
\centering
\includegraphics[width=1.0\textwidth]{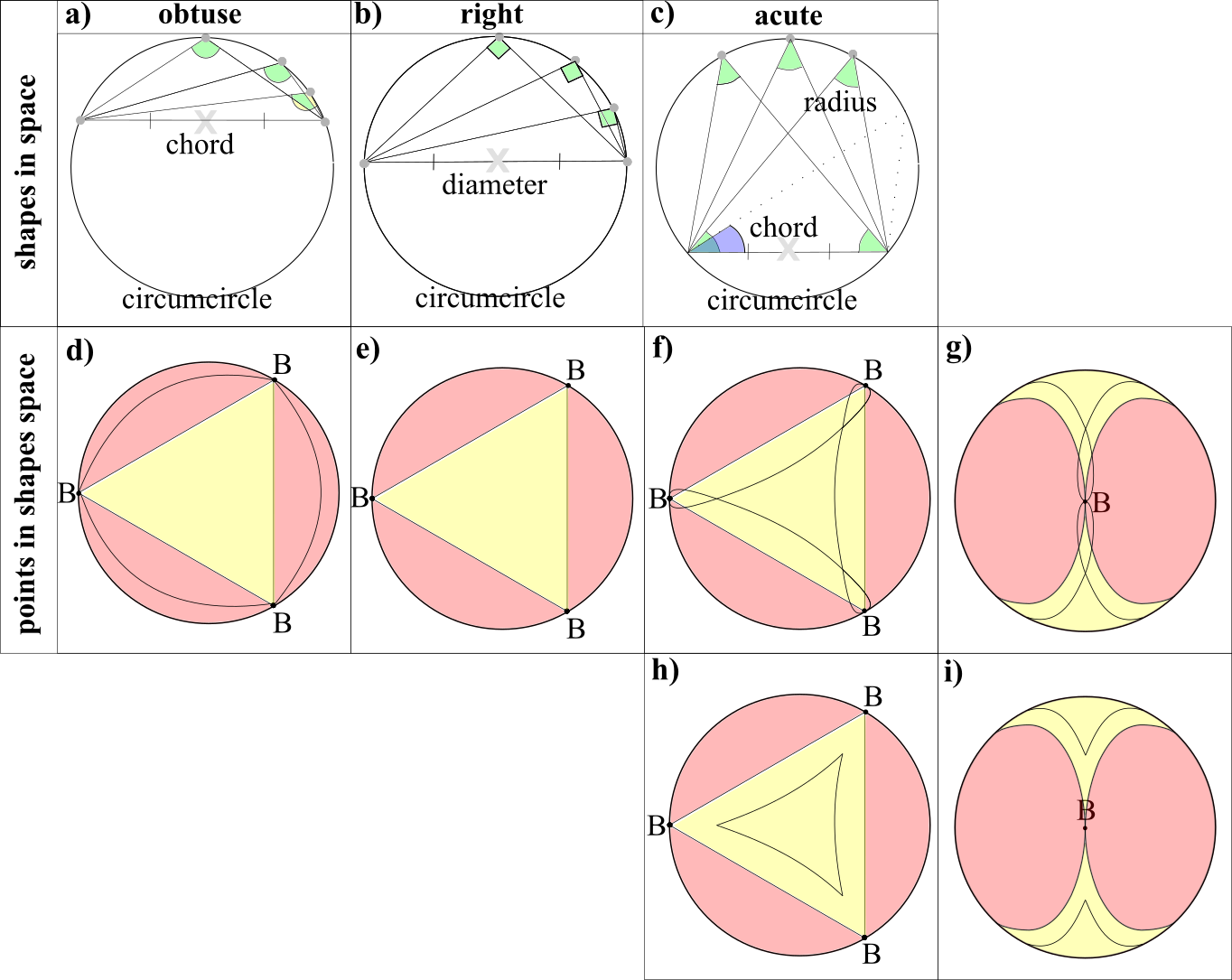}
\caption[Text der im Bilderverzeichnis auftaucht]{        \footnotesize{a) Constant obtuse-angle triangles in space and 
                                                                        d) the corresponding constant-obtuse-angle curves in shape space.

\mbox{ }

\ni b) Right triangles in space and 
    e) the corresponding kissing cap-circles of rightness in shape space, included for comparison.  

\mbox{ }

\ni c) Constant acute-angle triangles in space and 
    f) the corresponding constant-acute-angle curves in shape space. 
The acute angle shaded green is moreover maximal up to and including the isosceles triangles containing a second copy of it; 
beyond, in the dotted triangle depicted, the blue angle is maximal rather than the green one. 
This maximality is then reflected in shape space by discarding the part of the constant-angle curve beneath the intersection point, thus leaving a cusped simple curve as per h). 

\mbox{ }

\ni The construction used in a), b) and c) is based on the Chords Subtend Equal Angles Theorem of basic Euclidean Geometry.  
g) and i) provide a side-view of g) and i); this illustrates that clockwise and anticlockwise oriented constant-acute-angle curves `kiss' at B.  
However, as argued above, this is to be excluded due to its $\alpha$-undefined rather than constant-$\alpha_{\sm\sa\sx}$ status as a shape-in-space.} } 
\label{alpha-max} \end{figure}          }

\mbox{ }

\ni Let us next sketch and interpret these constant-$\alpha$ curves, 
as well as the consequent constant-$\alpha_{\sm\sa\sx}$ curves, for which there are three angles that could be maximal. 

\mbox{ } 

\ni{\bf Remark 1} On the one hand, for obtuse constant-$\alpha$,
\be 
\alpha \mbox{ } > \mbox{ } \frac{\pi}{2} \mbox{ } ,
\ee 
the shapes-in-space are as in Fig \ref{alpha-max}.a).
Each of the three labelling's version of the constant-angle curve is disjoint (Fig \ref{alpha-max}.b).
Establishing this disjointness includes using the following shapes-in-space argument. 
While these curves jointly terminate four branches at a time at the binary collisions B, 
these points are themselves to excluded from the constant-angle curves since angles are undefined there.  
Moreover, being obtuse constant-$\alpha$ always implies being a curve of constant-$\alpha_{\sm\sa\sx}$, 
as the constant angle sum of the triangle leaves no room for any larger angle.

\mbox{ }

\ni{\bf Remark 2} The equator of collinearity is included as a limiting case of maximal angle curve, as the value of that maximal angle tends to its own maximal value, $\pi$. 
This limiting case differs from the above description in that just two branches end at each binary collision B.  
Moreover, since angle-undefinedness continues to apply at B, these need to be excised, leaving one with the equator with three equally spaced-out punctures.  

\mbox{ }  

\ni{\bf Remark 3} On the other hand, for acute constant-$\alpha$, 
\be
\alpha \mbox{ } < \mbox{ } \frac{\pi}{2} \mbox{ } , 
\ee
the shapes-in-space are as in Fig \ref{alpha-max}.f-g).
Now each of the three labelling's version of the constant-angle curve intersects with the other two, 
due to the existence of isosceles triangles with a second angle $\alpha$ (labelled one of the remaining ways or the other).  
This is also indicated in Fig \ref{alpha-max}.f-g).  
Once again, the constant-$\alpha$ curves terminate at the binary collisions B, which are to be excluded due to angles being undefined there.  
For acute constant-$\alpha$, moreover, constant-$\alpha$ need not imply constant-$\alpha_{\sm\sa\sx}$, 
since, once one swings past the isosceles triangles with two $\alpha$ angles, an angle other than $\alpha$ becomes maximal (Fig \ref{alpha-max}.e) again).  
Thus one keeps the cusped curved triangle parts.  
Namely,          one in the Northern alias E or clockwise-oriented hemisphere, depicted in Fig \ref{alpha-max}.h), 
and an identical one in the Southern alias $\overline{\mE}$ or anticlockwise-oriented hemisphere;   
pieces of both are visible in Fig \ref{alpha-max}.i)'s side-view. 

\mbox{ }

\ni{\bf Remark 4} The above acute-$\alpha_{\sm\sa\sx}$ analysis is moreover implicitly for 
\be
\frac{\pi}{3} \mbox{ } < \mbox{ } \alpha_{\sm\sa\sx} 
              \mbox{ } < \mbox{ } \pi/2               \mbox{ } . 
\ee 
This is firstly because 
\be 
\alpha_{\sm\sa\sx} \mbox{ } < \mbox{ } \pi/3 \mbox{ } \mbox{ is impossible } , 
\ee 
from combining 
\be 
\alpha_{\sm\sa\sx} \mbox{ } \geq \mbox{ } \alpha, \beta, \gamma \mbox{ } \geq \mbox{ } 0
\label{max}
\ee 
and the angle-sum of a flat triangle, 
\be
\alpha + \beta + \gamma = \pi \mbox{ } . 
\ee 
Secondly, the saturated case  
\be 
\alpha_{\sm\sa\sx} \mbox{ } = \mbox{ } \frac{\pi}{3} 
\ee 
furthemore saturates the first inequality in (\ref{max}), forcing equilaterality. 
This gives a separate case, since now all labels are equivalent, 
and one has collapsed down from an $\alpha_{\sm\sa\sx}$ curve to two isolated $\alpha_{\sm\sa\sx}$ points: E and $\overline{\mE}$.   

\mbox{ }

\ni{\bf Remark 5} There is moreover a sense in which B belongs to the kissing caps of rightness but not to any of the other $\alpha_{\sm\sa\sx}$ curves.  
Namely, that B is the limit of the triangle in which two different right angles are realized, 
by which it is justified in lying on the intersection of two distinct curves of rightness.
On the other hand, there are no triangles containing two equal obtuse angles.
It is in this way that we justify our description of kissing cap-circles of rightness while, on the other hand, taking all obtuse-or-collinear 
$\alpha_{\sm\sa\sx}$ to exclude their B limit points.  

\mbox{ }

\ni{\bf Remark 6} Remark 2's study implies that the cusp points correspond to flat isosceles triangles $\mI^{\sF}$; 
since these form three meridians, the cusp points moreover align from contour to contour $\forall \, \alpha_{\sm\sa\sx} \, \in \left(\frac{\pi}{3}, \, \frac{\pi}{2}\right)$.  

\mbox{ }

\ni{\bf Remark 7} It is also straightforward to establish that the line of tall isosceles triangles $\mI^{\sT}$ 
is the sole location of nontrivial stationary points of the maximal angle flow. 
Since these also form three meridians, the stationary points also align from contour to contour, now $\forall \, \alpha_{\sm\sa\sx} \, \in \left(\frac{\pi}{3}, \, \pi\right)$.  
By `nontrivial stationary points', mean that we exclude everything being stationary for the collinear configurations' equator. 

\mbox{ }

{            \begin{figure}[!ht]
\centering
\includegraphics[width=1.0\textwidth]{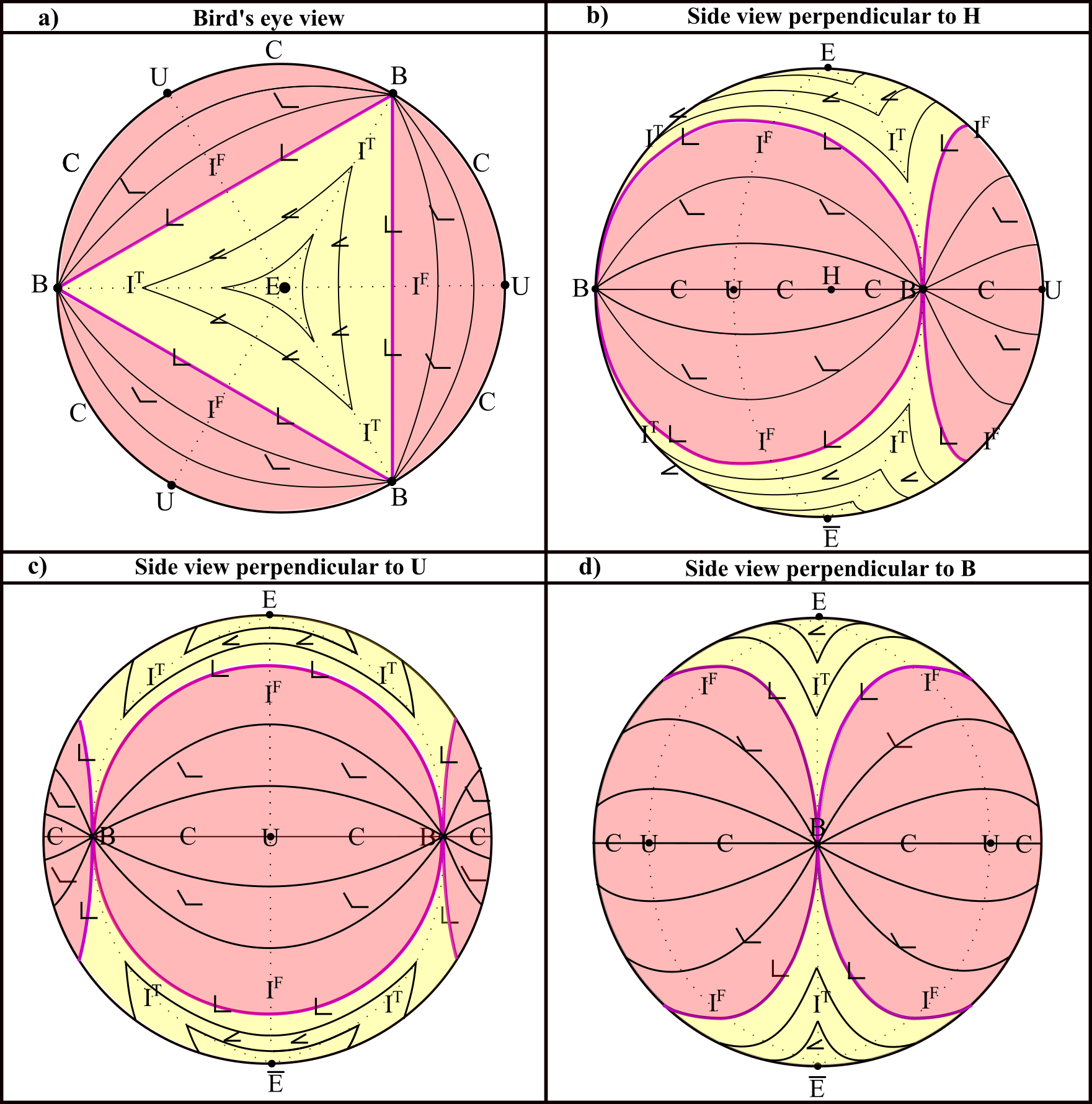}
\caption[Text der im Bilderverzeichnis auftaucht]{        \footnotesize{The maximal angle flow over the shape sphere of triangles as seen from various directions.
The kissing cap-circles of rightness separatrix is indicated in purple, with the two qualitatively different types of region separated by this shaded in red and yellow. 
The dotted lines indicate alignments of cusps and of stationary points.} }
\label{Maximal-Angle-Flow} \end{figure}          }

\vspace{10in}

\ni{\bf Remark 8} Remarks 1 to 5 establish that the kissing cap-circles of rightness constitute a {\it separatrix} for the maximal angle flow.   
I.e. a curve bounding between regions of qualitatively different behaviour.  
Moreover this separatrix in question is clearly not a simple curve, possessing, rather, three self-intersections at the B-points.  

\mbox{ } 

\ni{\bf Remark 9} The qualitatively different behaviours separated in the current problem are as follows. 

\mbox{ }

\ni a) The two acute regions, corresponding to the maximal angles 
\be 
\frac{\pi}{3} \mbox{ } < \mbox{ } \alpha_{\sm\sa\sx} \mbox{ } < \mbox{ } \frac{\pi}{2} \mbox{ } , 
\ee 
for which the maximal angle flow is cusped.

\mbox{ }

\ni b) The three obtuse regions, corresponding to the maximal angles
\be 
\frac{\pi}{2} \mbox{ } < \mbox{ } \alpha_{\sm\sa\sx} \mbox{ } \leq \mbox{ } \pi \mbox{ } . 
\ee
for which limit points of almost all maximal-angle curves are not part of those curves, the exception being the bounding case of the equator of collinearity.  

\mbox{ }

\ni{\bf Remark 10} So we have found a geometrically-meaningful example of a separatrix separating between two different kinds of analytic pathology of curves: 
between curves with nondifferentiable cusps and curves whose limiting points are not defined as part of the curve. 
This example is likely to be of interest in both Differential Geometry and Dynamical Systems, 
as well as in Flat Geometry due to its correspondence to interesting new elementarily phraseable (while far from necessarily elementarily solvable) propositions about triangles.

\mbox{ }

\ni{\bf Remark 11} The central point of each acute region  -- the equilateral triangles of each orientation: E and $\overline{\mE}$ -- is qualitatively distinct. 
These form $\alpha_{\sm\sa\sx}$ orbits by themselves $\left(\mbox{for the value $\frac{\pi}{3}$}\right)$, whereas all other $\alpha_{\sm\sa\sx}$ orbits are (collections of) curves. 
As critical points, the equilateral triangles  E and $\overline{\mE}$ have moreover the status of {\it centres} in the maximal-angle flow. 
They are moreover nonlinear, indeed non-differentiable, centres due to their flowlines containing three cusps each.  

\mbox{ }

\ni{\bf Remark 12} On the other hand, the B points: binary collisions at which kissing occurs, are considerably more complicated in their further aspect as critical points 
in the maximal angle flow, as described in Fig \ref{B-as-Critical-Point}. 
%
{            \begin{figure}[!ht]
\centering
\includegraphics[width=0.32\textwidth]{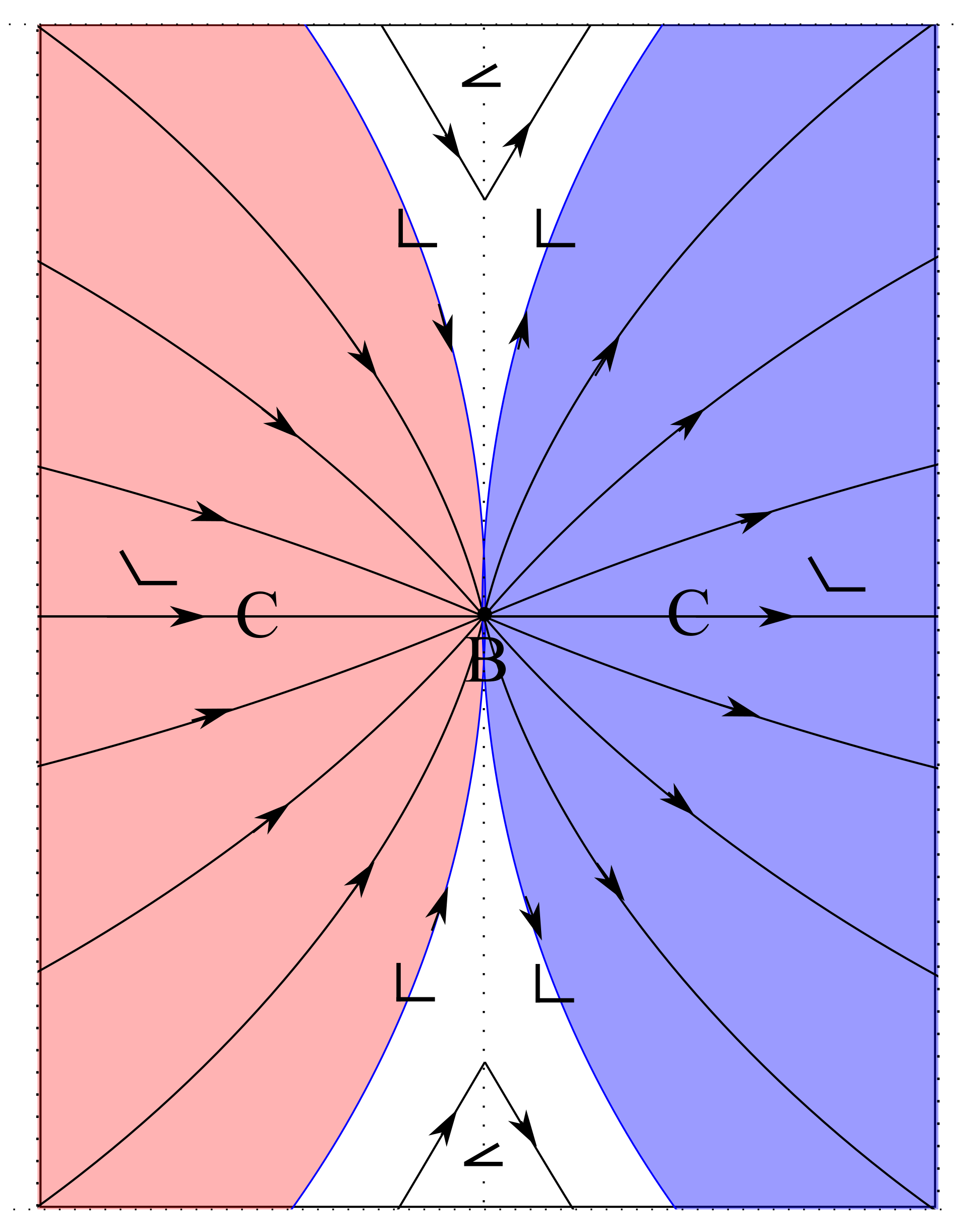}
\caption[Text der im Bilderverzeichnis auftaucht]{        \footnotesize{Detail of B as a critical point. 
Note it has an attracting node wedge subregion (red) and a repelling node wedge subregion (blue), with `cusped saddle' wedge subregions (white) between the two.} }
\label{B-as-Critical-Point} \end{figure}          }

\mbox{ }

\ni{\bf Remark 13} All in all, the maximal angle flow has five critical points: three B's (one per clustering choice) and two E's (one per choice of orientation).            

\mbox{ }

\ni{\bf Remark 14} Finally, this flow is reflection-symmetric about the equator.
It is also $S_3 = D_3$ symmetric as regards exchanging the B's. 
$S_3$ refers here to the permutation group acting on the vertex labels, 
whereas $D_3$ is the dihaedral group of symmetries of the bird's-eye-view triangle in shape space formed by the cap-circles of rightness.  
Thus the overall symmetry group is $S_3 \times \mathbb{Z}_2$, corresponding identically to the symmetry group of the shape sphere itself.

\vspace{10in}

\section{Probability that a triangle has nontrivially-located Fermat point}

\ni{\bf Corollary 3} For each cluster i) 
\be
Area \mbox{ } = \mbox{ } 
2 \, \left| \int_{{\cal X} = 0}^{1/2} \d {\cal X} \, \mbox{arcsin}  \left(  \frac{    1 - 2 \, {\cal X}    }{    \sqrt{3} \, k \, \sqrt{   1 - {\cal X}^2  }   }  \right) \right| 
\mbox{ } . 
\ee
\ni{\underline{Proof}} Using a multiplicative factor of 2 for the upper and lower hemispheres as related by the reflection symmetry about the equator, 
\be 
Area \mbox{ } = \mbox{ } 2 \, \left|  \int_{\Theta = 0}^{\pi/3} \d \Theta \, \mbox{sin} \, \Theta \int_{\Phi = 0}^{\Phi(\Theta)} \d \Phi  \right| 
     \mbox{ } = \mbox{ } 2 \, \left|  \int_{{\cal X} = 0}^{1/2} \d {\cal X}                       \int_{\Phi = 0}^{\Phi(\cal X)} \d \Phi  \right| \mbox{ } , 
\ee 
where $\Phi({\Theta})$ is given by (\ref{Phi-eq-2}) and consequently $\Phi({\cal X})$ is given by (\ref{Phi-eq-3}). $\Box$

\mbox{ }
 
\ni{\bf Remark 15} The underlying indefinite integral here is in fact analytically tractable \cite{Maple}, 
in terms of complete elliptic integrals of the first, second and third kind \cite{Ellip-Fns}. 
The area calculation itself moreover involves a definite integral which is numerically tractable \cite{Maple}. 

\mbox{ }

\ni{\bf Remark 16} We give the following as a particularly interesting concrete example.  

\mbox{ }

\ni{\bf Corollary 4} for the Fermat value $\alpha_{\sm\sa\sx} = 2 \, \pi/3$, 
\be
Area \mbox{ } = \mbox{ }  2 \left| \int_{{\cal X} = 0}^{1/2} \d {\cal X} \, 
\mbox{arcsin}  \left(  \frac{    2 \, {\cal X} - 1    }{   \, \sqrt{   1 - {\cal X}^2  }   }  \right) \right| \mbox{ } . 
\ee
This can moreover be evaluated using Maple \cite{Maple}. 

\mbox{ }

\ni{\bf Remark 17} Furthermore, once again normalizing by the total area of the sphere, 
\be
\mbox{Prob($\alpha$-obtuse)} \mbox{ } = \mbox{ }  \frac{3 \, Area}{4 \, \pi} \mbox{ } . 
\ee
Thus 
\be 
\mbox{Prob(Fermat-obtuse)} = 0.1394 \mbox{ } . 
\ee
Complementarily, 
\be
\mbox{Prob($\alpha$-acute)} =  1 - \mbox{Prob($\alpha$-obtuse)} = \frac{4 \, \pi - 3 \, Area}{4 \, \pi} \mbox{ } , 
\ee 
So 
\be 
\mbox{Prob(Fermat-acute)} = 0.8606 \mbox{ } .  
\ee
It is this last quantity which is of particular interest due to its correspondence to the triangles whose Fermat points are nontrivially positioned.  

\vspace{10in}

\section{Conclusion}

\ni Lewis Carroll's longstanding and conceptually interesting pillow problem of what is the probability that a triangle is obtuse 
has recently been solved \cite{MIT, III, A-Pillow} using Kendall's Shape Theory \cite{Kendall84, Kendall89, Kendall}.
In the current paper, we have presented a substantial extension of this problem and its shape-theoretic solution; 
\cite{A-Pillow, IV} give distinct extensions, indicating this to be a {\sl fertile} source of new research in Geometry.  
The current paper's extension follows from viewing right-angled triangles as triangles whose maximal angle is right, 
and then passing to consider all the other possible values of maximal angle, $\alpha_{\sm\sa\sx} = [\frac{\pi}{3}, \, \pi]$.
This gives the {\it maximal angle flow} over Kendall's shape sphere.  

\mbox{ }

\ni We find that the three kissing cap-circles of rightness 
-- that so readily provide the shape-theoretic resolution of Lewis Carroll's pillow problem of the probability that a triangle is obtuse -- 
furthermore play the role of separatrix in the maximal angle flow. 
The extent to which right angles' significance in the geometry for triangles is underlied by this triangleland shape space separatrix result remains to be determined.  
The qualitatively different behaviour in each of the regions separated thus is moreover very interesting. 

\mbox{ }

\ni On the one hand, in the acute regions, the flowlines have three cusps. 

\mbox{ }

\ni On the other hand, in the obtuse regions, the flowlines have limit points which do not belong to that flowline.  

\mbox{ }

\ni Thus the kissing cap-circles of rightness are a separatrix between flow regions which each exhibit a distinct analytic pathology. 

\mbox{ }

\ni In the maximal angle flow, the equilateral triangle points E, $\overline{\mE}$ have moreover the status of centres (nonlinear ones with 3 cusps on the circulating curves).  
In contrast, the B points -- binary collision shapes and the kissing points of the cap-circles of rightness -- 
are more complicated critical points, as per Fig \ref{B-as-Critical-Point}.

\mbox{ }

\ni We apply our new knowledge of the maximal angle flow over the shape sphere to find the probability that a triangle has nontrivially-located Fermat point. 
This corresponds to the critical angle $\alpha_{\sm\sa\sx} = \frac{2 \, \pi}{3}$. 
This is critical in the sense that if the angle is smaller than this (`Fermat acute'), 
the Fermat point lies in the interior of the triangle, but if it is larger (`Fermat obtuse'), the Fermat point just lies on the obtuse vertex.
We obtain $\mbox{Prob(Fermat-acute} = 0.8606$.  

\mbox{ }

\ni Let us end by pointing to various extensions of this work.

\mbox{ }

\ni Extension 1) Our flow and area-ratio probability method can be applied to find further maximal critical angle dependent Flat Geometry triangles' properties' probabilities 
of occurring in the shape-theoretic sense of `random triangles'.

\mbox{ }

\ni Extension 2) The current paper's maximal angle flow over the shape sphere's above-described features likely render the current paper of interest as an example 
in Differential Geometry and Dynamical Systems. 

\mbox{ }

\ni Extension 3) \cite{MIT, III, A-Pillow} and the current paper's use of Shape Theory to answer questions about random triangles moreover extends to 
its use in answering questions about random quadrilaterals and polygons using the Differential Geometry of the $N$-a-gonland $\mathbb{CP}^{N - 2}$ \cite{Kendall, QuadI, IV}. 
This program uses Differential Geometry on Kendall's shape spaces to pose and/or solve probabilistic questions 
involving polygons by considering these to live on Kendall's natural shape space with uniform measure thereupon.
Fermat's own problem has moreover a large number of generalizations as well -- geometric means -- which overlap well with the theory of polygons. 

\mbox{ }

\ni Extension 4) The current paper's notion of similarity Shape Theorey moreover itself extends to a wide variety of further notions of shape and shape space, 
such as affine, projective, conformal, M\"{o}bius... \cite{GT09, AMech, ABook, PE16, KKH16}.   
Many of Extension 3)'s considerations recur in this wider arena.  

\mbox{ }

\ni{\bf Acknowledgments} I thank Chris Isham and Don Page for previous discussions, 
and Reza Tavakol, Malcolm MacCallum, Enrique Alvarez and Jeremy Butterfield for support with my career.



\begin{thebibliography}{99}

\footnotesize

\bibitem{17th}                See R.A. Johnson, {\it Modern Geometry: An Elementary Treatise on the Geometry of the Triangle and the Circle} (Houghton Mifflin, Boston 1929) 
                              pp. 221-222 for discussion of the original 17th Century contributions of Fermat and Torricelli in this regard.  

\bibitem{Pillow}              C.L. Dodgson (alias Lewis Carroll), {\it Curiosa Mathematica: Pillow-Problems, thought out during Sleepless Nights} (Macmillan, London 1893).							  

\bibitem{Ellip-Fns}           M. Abramowitz and I.A. Stegun, {\it Handbook of Mathematical Functions with Formulas, Graphs, and Mathematical Tables} (Dover, New York 1972). 
 						  
\bibitem{Kendall84}           D.G. Kendall, ``Shape Manifolds, Procrustean Metrics and Complex Projective Spaces", Bull. Lond. Math. Soc. {\bf 16} 81 (1984). 

\bibitem{Kendall89}           D.G. Kendall, ``A Survey of the Statistical Theory of Shape", Statistical Science {\bf 4} 87 (1989).

\bibitem{Marchal}             C. Marchal, {\it Celestial Mechanics} (Elsevier, Tokyo 1990).

\bibitem{Guy}                 R. Guy, ``There Are Three Times as Many Obtuse-Angled Triangles as There Are Acute-Angled Ones",
                              Mathematics Magazine, {\bf 66} 175 (1993).   

\bibitem{generalizations}     G. Wesolowsky, {\it The Weber problem: History and Perspective}, Location Science {\bf 1} 5 (1993).
							  
\bibitem{Portnoy}             S. Portnoy, ``A Lewis Carroll Pillow Problem: Probability of an Obtuse Triangle", Statist. Sci. {\bf 9} 279 (1994).  
			
\bibitem{Small}               C.G.S. Small, {\it The Statistical Theory of Shape} (Springer, New York, 1996).  
							  			  						  
\bibitem{Kendall}             D.G. Kendall, D. Barden, T.K. Carne and H. Le, {\it Shape and Shape Theory} (Wiley, Chichester 1999).  
												   				
\bibitem{Montgomery}          R. Montgomery, ``Fitting Hyperbolic Pants to a 3-Body Problem", Ergod. Th. Dynam. Sys. {\bf 25} 921 (2005), math/0405014.

\bibitem{FORD}                E. Anderson, ``Foundations of Relational Particle Dynamics", Class. Quant. Grav. {\bf 25} 025003 (2008), arXiv:0706.3934.
																
\bibitem{+Tri}                E. Anderson, ``Shape Space Methods for Quantum Cosmological Triangleland", Gen. Rel. Grav. {\bf 43} 1529 (2011), arXiv:0909.2439.  

\bibitem{GT09}                D. Groisser, and H.D. Tagare, ``On the Topology and Geometry of Spaces of Affine Shapes", 
                              Journal of Mathematical Imaging and Vision {\bf 34} 222 (2009).  

\bibitem{Frankel}             T. Frankel, {\it The Geometry of Physics: An Introduction} (Cambridge University Press, Cambridge 2011).  
							  
\bibitem{FileR}               E. Anderson, ``The Problem of Time and Quantum Cosmology in the Relational Particle Mechanics Arena", arXiv:1111.1472.  

\bibitem{QuadI}               E. Anderson, ``Relational Quadrilateralland. I. The Classical Theory", Int. J. Mod. Phys. {\bf D23} 1450014 (2014), arXiv:1202.4186.

\bibitem{Bhatta}              A. Bhattacharya and R. Bhattacharya, {\it Nonparametric Statistics on Manifolds with Applications to Shape Spaces} 
                             (Cambridge University Press, Cambridge 2012).
		  
\bibitem{QuadII}              E. Anderson and S.A.R. Kneller, ``Relational Quadrilateralland. II. The Quantum Theory", Int. J. Mod. Phys. {\bf D23} 1450052 (2014), arXiv:1303.5645.

\bibitem{MIT}                 A. Edelman and G. Strang, ``Random Triangle Theory with Geometry and Applications", 
                              Foundations of Computational Mathematics (2015), arXiv:1501.03053.   
							  				  
\bibitem{AMech}               E. Anderson, ``Six New Mechanics corresponding to further Shape Theories", Int. J. Mod. Phys. {\bf D 25} 1650044 (2016), arXiv:1505.00488. 

\bibitem{DM16}                I.L. Dryden, K.V. Mardia, {\it Statistical Shape Analysis: With Applications in R}, 2nd Edition (Wiley, Chichester 2016).  
	
\bibitem{PE16}                V. Patrangenaru and L. Ellingson ``Nonparametric Statistics on Manifolds and their Applications to Object Data Analysis" 
                             (Taylor and Francis, Boca Raton, Florida 2016).  
			 
\bibitem{KKH16}               F. Kelma, J.T. Kent and T. Hotz, ``On the Topology of Projective Shape Spaces", arXiv:1602.04330. 
		
\bibitem{ABook}               E. Anderson, {\it Problem of Time. Quantum Mechanics versus General Relativity}, (Springer International 2017), Fundam.Theor.Phys. 190 (2017) pp.-;  
                              free access to its extensive Appendices can be found at https://link.springer.com/content/pdf/bbm$\%$3A978-3-319-58848-3$\%$2F1.pdf .
			
\bibitem{I}                   E. Anderson, ``The Smallest Shape Spaces. I. Shape Theory Posed, with Example of 3 Points on the Line", arXiv:1711.10054.

\bibitem{II}                  E. Anderson, ``The Smallest Shape Spaces. II. 
                              4 Points on a Line Suffices for a Complex Background-Independent Theory of Inhomogeneity", arXiv:1711.10073.

\bibitem{III}                 E. Anderson, ``The Smallest Shape Spaces. III. Triangles in the Plane and in 3-$d$", arXiv:1711.10115.
	
\bibitem{A-Pillow}            E. Anderson, ``Alice in Triangleland: Lewis Carroll's Pillow Problem and Variants Solved on Shape Space of Triangles", arXiv:1711.11492.

\bibitem{2-Herons}            E. Anderson, ``Two New Perspectives on Heron's Formula", arXiv:1712.01441.  
							  
\bibitem{Ineq}                E. Anderson, ``Shape (In)dependent Inequalities for Triangleland's Jacobi and Democratic-Linear Ellipticity Quantitities", arXiv:1712.04090.
							   
\bibitem{Maple}               This was calculated using Maple 17.

\bibitem{A-Monopoles}         E. Anderson, ``Monopoles of Fifteen Types in 3-Body Problems", forthcoming 2017.  

\bibitem{A-Perimeter}         E. Anderson, ``Geometrically Significant Shape Quantities for Triangles Extremized over the Triangleland Shape Sphere", forthcoming 2018.  

\bibitem{IV}                  E. Anderson, ``The Smallest Shape Spaces. IV. Quadrilaterals in the Plane", forthcoming 2018.
 
\bibitem{Forthcoming}         E. Anderson, forthcoming.
 
\end{thebibliography}
\end{document}